\documentclass{amsart}
\usepackage{color}
\usepackage[colorlinks]{hyperref}
\usepackage{tikz}
\usetikzlibrary{arrows}

\newtheorem{thm}{Theorem}[section]
\newtheorem{cor}[thm]{Corollary}
\newtheorem{lem}[thm]{Lemma}

\theoremstyle{definition}
\newtheorem{defn}[thm]{Definition}
\theoremstyle{remark}
\newtheorem{rem}[thm]{Remark}
\numberwithin{equation}{section}

\begin{document}
\title{ Equivalencies between beta-Shifts and $S$-gap Shifts}
\author{D. Ahmadi Dastjerdi and S. Jangjooye Shaldehi }%
\address{Department of Pure Mathematics, Faculty of Mathematical Sciences ,The University of Guilan}
\email{dahmadi1387@gmail.com, sjangjoo90@gmail.com}%


\begin{abstract}
Let $ X_{\beta}$ be a sofic $ \beta $-shift for $ \beta \in (1,\,2] $. We show that there is an $ S $-gap shift $ X(S) $ such that $ X_{\beta} $ and $ X(S) $ are right-resolving almost conjugate.
Conversely,  a condition on $ S \subseteq \mathbb N\cup \{0\} $ is given such that  for this $S$, there is a $ \beta $ such that  $ X(S) $ and $ X_{\beta} $ have the same equivalency.
We  show that if $X_\beta$ is SFT, then there is an $S$-gap shift  conjugate to this $X_\beta$; however, if $X_\beta$ is not SFT, then no $S$-gap shift  is conjugate to  $X_\beta$. Also
we will investigate the existence of these sort of equivalencies for  non-sofics.
\end{abstract}

\maketitle
\section*{Introduction}
Two important classes of symbolic dynamics are $ \beta $-shifts and $S$-gap shifts.
Both are rich families and highly chaotic with application in coding theory, number theory and a source of examples for symbolic dynamics.
There has been some independent studies of these classes. See \cite{CT, P, T1, T2} for $ \beta $-shifts and \cite{A1, CT, Jung} for $S$-gap shifts. There are some common properties among these two classes. For instance, both of them are coded systems having positive entropies for $ |S|\geq 2 $ and every subshift factor of them is intrinsically ergodic \cite{CT}. There are disparities as well: $ \beta $-shifts are all mixings, though this is not true for $S$-gap shifts \cite{Jung}; and $S$-gap shifts are synchronized, a property which does not hold for all $ \beta $-shifts. Even among sofic $ \beta $-shifts and $S$-gap shifts, which are of our most interest here, there are some major differences. An important class of sofic $S$-gap shifts are almost-finite-type (AFT) \cite{A1}, but no $ \beta $-shift is AFT \cite{T2}.

Here we let $\beta\in(1,\,2]$ and search for an $S$-gap shift $X(S)$ which has some sort of equivalencies with our $ \beta $-shift denoted by $X_\beta$. A well known equivalency between dynamical systems is conjugacy and in transitive systems, almost conjugacy which is virtually a conjugacy between transitive points has been much considered. However, in coding theory the equivalencies which are  right-resolving - called \emph{deterministic} in computer science - are more natural and applicable. Our systems can have application in coding theory and so we investigate the equivalency of this sort among them.

Here, we summarize our results. Let $\beta\in (1,\,2]$ and let $X_\beta$ be the corresponding $\beta$-shift.
We will associate to $X_\beta$ a unique $S$-gap shift denoted by ASS$(X_\beta)$ and 
 will show that  when $X_\beta$ is
sofic, then  $X_{\beta} $ and $ X(S) =$ASS$(X_\beta)$ are
 right-resolving almost conjugate and when $X_\beta$ is SFT, they are conjugate as well (Theorem  \ref{finite equivalent}).  On the other side,
for a given  $S$-gap shift, it does not necessarily exist a $ \beta $-shift holding the same equivalencies with X(S). However, we give
a necessary and sufficient condition on $S$
   to have this true (Theorem \ref{inverse}). 
   
   In section 4, we extend the results obtained for sofics to  non-sofics. For instance,  $ X_{\beta} $ and ASS$(X_{\beta}) $ have a common extension 
which is an almost Markov synchronized system whose maps are entropy-preserving (Corollary \ref{almost Markov}). Additionally, when $ X_{\beta} $ is synchronized, $ X_{\beta} $ and ASS$(X_{\beta}) $ have a common synchronized 1-1 a.e. extension with right-resolving legs (Theorem \ref{common synchronized}). 

Theorem \ref{zeta} gives the the zeta function of $X_\beta$ in terms of the zeta function of $\mbox{ASS}(X_\beta)$ and
Theorem \ref{cantor set} states that $\{\mbox{ASS}(X_\beta):\, \beta\in (1,\,2]\}$ is a Cantor set in the set of all $S$-gap shifts.

\section{Background and Notations}
The notations has been taken from \cite{LM} and the proofs of the
claims in this section can be found there. Let $\mathcal A$ be an
alphabet, that is a nonempty finite set. The full $\mathcal
A$-shift
 denoted by ${\mathcal A}^{\mathbb Z}$, is the collection of all bi-infinite sequences of symbols from ${\mathcal A}$.
A block (or word) over ${\mathcal A}$ is a finite sequence of
symbols from ${\mathcal A}$. The shift function $\sigma$ on the full
shift ${\mathcal A}^{\mathbb Z}$ maps a point $x$ to the point
$y=\sigma(x)$ whose $i$th coordinate is $y_{i}=x_{i+1}$.

 Let ${\mathcal B}_{n}(X)$ denote the set of all admissible $n$-blocks.
 The \emph{Language} of $X$ is the collection ${\mathcal B}(X)=\bigcup_{n=0}^{\infty}{\mathcal B}_{n}(X)$.
A word $v \in {\mathcal B}(X)$ is \emph{synchronizing} if whenever $uv$ and $vw$ are in ${\mathcal B}(X)$,
we have $uvw \in {\mathcal B}(X)$.

Let ${\mathcal A}$ and ${\mathcal D}$ be alphabets and $X$ a
shift space over ${\mathcal A}$. Fix integers $m$ and $n$ with $-m
\leq n$. Define the \emph{$(m+n+1)$-block map} $\Phi: {\mathcal
B}_{m+n+1}(X) \rightarrow {\mathcal D}$ by
\begin{equation}\label{sliding} 
y_{i}=\Phi(x_{i-m}x_{i-m+1}...x_{i+n})=\Phi(x_{[i-m,i+n]})
\end{equation}
where $y_{i}$ is a symbol in ${\mathcal D}$. The map $\phi: X
\rightarrow {\mathcal D}^{\mathbb Z}$ defined by $y=\phi(x)$ with
$y_{i}$ given by \eqref{sliding} is called the \emph{sliding block code}
with \emph{memory} $m$ and \emph{anticipation} $n$ induced by
$\Phi$. An onto sliding block code $\phi: X \rightarrow Y$ is called a
factor code. In this case, we say that $Y$ is a
factor of $X$. The map $\phi$ is a \emph{conjugacy}, if it is invertible.

 An \emph{edge shift}, denoted by $X_{G}$,  is a shift space which consists of all
  bi-infinite walks in a directed graph $G$.
A \emph{labeled graph} ${\mathcal G}$ is a pair $(G,{\mathcal L})$ where
$G$ is a graph with edge set $\mathcal E$, vertex set $\mathcal V$ and the labeling
${\mathcal L}: {\mathcal E} \rightarrow {\mathcal A}$. Each $ e \in \mathcal E $ starts at a vertex denoted by $i(e) \in \mathcal V   $ and terminates at a vertex $t(e) \in \mathcal V    $.

When the set of forbidden words is finite, the space is called \emph{subshift of finite type} (SFT).
A
\emph{sofic shift} $X_{{\mathcal G}}$ is the set of sequences
obtained by reading the labels of walks on $G$,
$$
 X_{{\mathcal G}}=\{{\mathcal L}_{\infty}(\xi): \xi \in X_{G}\}={\mathcal L}_{\infty}(X_{G}).
$$
We say ${\mathcal G}$ is a \emph{presentation} of $X_{{\mathcal
G}}$. 

A labeled graph ${\mathcal G}=(G,{\mathcal L})$ is
\emph{right-resolving} if for each vertex $I$ of $G$ the edges
starting at $I$ carry different labels. A \emph{minimal
right-resolving presentation} of a sofic shift $X$ is a
right-resolving presentation of $X$
 having the fewer vertices among all right-resolving presentations of $X$. Any two minimal right-resolving
presentations of an irreducible sofic shift must be isomorphic as
labeled graphs \cite[Theorem 3.3.18]{LM}. So we can speak of ``the" minimal right-resolving presentation
of an irreducible sofic shift $X$ which we call it the \emph{Fischer cover} of $ X $.

Let $w \in {\mathcal B}(X)$. The
\emph{follower set} $F(w)$ of $w$ is defined by $F(w)=\{v \in
{\mathcal B}(X): wv \in  {\mathcal B}(X)\}.$ A shift space $X$ is
sofic if and only if it has a finite number of follower sets
\cite[Theorem 3.2.10 ]{LM}. In this case, we have a labeled graph
${\mathcal G}=(G,{\mathcal L})$ called the \emph{follower set
graph} of $X$. The vertices of $G$ are the follower sets and if
 $wa \in {\mathcal B}(X)$, then draw an edge labeled $a$ from $F(w)$ to $F(wa)$. If $wa \not \in {\mathcal B}(X)$ then do nothing.

 Let $ \phi=\Phi_{\infty}:X\rightarrow Y $ be
a 1-block code. Then $ \phi $ is \emph{right-resolving} if whenever $ ab $ and $ ac $  are 2-blocks
in $ X $ with $ \Phi(b)=\Phi(c) $, then $ b=c $. 

Let $  G$ and $  H$ be graphs. A \emph{graph homomorphism} from $  G$
to $  H$ consists of a pair of maps $\partial \Phi:\mathcal{V}(G)\rightarrow \mathcal{V}(H)$  and $\Phi:\mathcal{E}(G)\rightarrow \mathcal{E}(H)$
such that $ \partial\Phi(i(e))=i(\Phi(e))  $ and
$ \partial\Phi(t(e))=t(\Phi(e)) $ for all $ e \in \mathcal E(G) $.
A graph homomorphism is a \emph{graph isomorphism} if both $ \partial \Phi $ and $ \Phi $ are one-to-one
and onto. Two graphs $  G$ and $  H$ are
graph isomorphic (written $  G \cong H$) if there is a graph isomorphism between
them. Let  $\mathcal E_{I}(G)  $ be the set of all the edges in $\mathcal{E}(G)$ starting from $I\in\mathcal{V}(G)$.
A graph homomorphism $ \Phi: G\rightarrow H $ maps $\mathcal E_{I}(G)  $ into $ \mathcal E_{\partial\Phi(I)}(H) $ for each
vertex $  I$ of $  G$. Thus $ \phi=\Phi_{\infty} $ is right-resolving if and only if for every
vertex $  I$ of $  G$ the restriction $\Phi_{I}  $ of $  \Phi$ to $\mathcal E_{I}(G)  $ is one-to-one. If $  G$ and $  H$ are irreducible and $ \phi $ is a right-resolving code from $ X_{G} $ onto $ X_{H} $, then each $\Phi_{I}  $ must be a bijection. Thus for each vertex $  I$ of $G  $ and every edge $f \in \mathcal E_{\partial\Phi(I)}(H) $,
there exists a unique ``lifted" edge $e \in \mathcal E_{I}(G)  $ such that $ \Phi(e)=f $. This lifting property inductively extends to paths: for every vertex $  I$ of $G  $ and every path $  w$ in $  H$ starting at $ \partial\Phi(I) $, there is a unique
path $  \pi$ in $  G$ starting at $  I$ such that $ \Phi(\pi)=w $.

Points $ x$ and $ x' $ in $ X $ are left-asymptotic
if there is an integer $ N $ for which $ x_{(-\infty,\,N]}=x'_{(-\infty,\,N]} $. A sliding block
code $ \phi: X\rightarrow Y $ is \emph{right-closing} if whenever $ x,\, x' $ are left-asymptotic and
$ \phi(x)=\phi(x') $, then $ x=x' $. 

 The entropy of a shift space $X$ is defined by $h(X)=\lim_{n
\rightarrow \infty}(1/n)\log|{\mathcal B}_{n}(X)|$.

\section{General Properties of $S$-gap Shifts and $ \beta $-shifts}\label{properties}
\subsection{$S$-gap shifts}\label{subs:S-gap}
To define an $S$-gap shift $X(S)$, fix $S=\{s_i\in\mathbb N\cup\{0\}:0\leq s_i<s_{i+1},\, i\in\mathbb{N} \cup \{0\}\}$. Define
$X(S)$ to be the set of all binary sequences for which $1$'s
occur infinitely often in each direction and such that the number
of $0$'s between successive occurrences of a $1$ is in
$S$. When $S$ is infinite, we need to allow points that begin or
end with an infinite string of $0$'s. Note that  $X(S)$ and $X(S^{'})$ are conjugate if and only if one of the $S$ and $S^{'}$ is $\{0,\,n\}$ and the other $\{n,\,n+1,\,n+2,\ldots\}$ for some $n \in \mathbb N$ \cite[Theorem 4.1]{A1}. So we consider $ X(S) $ up to conjugacy and by convention $\{0,n\}$ is chosen. Now let $d_{0}=s_{0}$ and $\Delta(S)=\{d_{n}\}_{n}$ where
$d_{n}=s_{n}-s_{n-1}$. Then an $S$-gap shift is subshift of finite type (SFT) if and only if $S$ is
finite or cofinite, is almost-finite-type (AFT) if and only if $\Delta(S)$ is
eventually constant and is sofic if and only if $\Delta(S)$ is
eventually periodic \cite{A1}. Therefore, for sofic $S$-gap shifts  we set
\begin{equation}\label{Delta}
\Delta(S)=\{d_{0},\,d_{1},\ldots,\,d_{k-1},\overline
{g_{0},g_{1},\ldots,g_{l-1}}\},\qquad g=\sum_{i=0}^{l-1} g_i
\end{equation}
where $g_{j}=s_{k+j}-s_{k+j-1},0\leq j \leq l-1$. Also $ k $ and $ l $ are the least integers such that \eqref{Delta} holds.

The Fischer cover of any irreducible sofic shift as well as $S$-gap shifts is the labeled subgraph of the follower set graph consists of the finite set of follower sets of synchronizing words as its vertices. For an $S$-gap shift this set is
\begin{equation}\label{follower}
 \{F(1),\,F(10),\ldots,\,F(10^{n(S)})\},
\end{equation}
where $n(S)=\max S$
for $|S|<\infty$.
 If $|S|=\infty$, then $ n(S) $ will be defined as follow.
\begin{enumerate}
\item
For $k=1$ and $g_{l-1}> s_{0}$,
\begin{enumerate}
\item
if $g_{l-1}=s_{0}+1$,  then $ F(10^{s_{l-1}+1})=F(1) $ and $ n(S)=s_{l-1} $.
\item
if $g_{l-1}>s_{0}+1$, then $ F(10^{g})=F(1) $ and $ n(S)=g-1 $.
\end{enumerate}
\item
For $ k \neq 1 $, if $g_{l-1}>d_{k-1}$, then $ F(10^{g+s_{k-2}+1})=F(10^{s_{k-2}+1}) $ and $ n(S)=g+s_{k-2} $.
\item
For $ k \in \mathbb N $, if $g_{l-1}\leq d_{k-1}$, then $ F(10^{s_{k+l-2}+1})=F(10^{s_{k-1}-g_{l-1}+1}) $ and $ n(S)=s_{k+l-2} $.
\end{enumerate}
To have a  view  about the Fischer cover of an $S$-gap shift,
 we line up vertices in \eqref{follower} horizontally starting from $F(1)$ on the left followed by $F(10)$ then by $F(10^2)$ and at last ending at
 $ F(10^{n(S)}) $ as the far right vertex. In all cases, label 0 the edge starting from $ F(10^{i}) $ and terminating at $ F(10^{i+1}) $, $ 0\leq i \leq n(S)-1 $; also, label $  1$ all edges from $F(10^s)$ to $F(1)$ for $s\in S$ and $s<n(S)$. 

So the only remaining edges to be taken care of are those starting at $F(10^{n(S)})$. In (1a), there are two edges from $ F(10^{n(S)}) $ to $ F(1) $; label one $  0$ and the other $  1$. 
In (1b), there is only one edge from $ F(10^{n(S)}) $ to $ F(1) $ which is  labeled $ 0 $.
In case (2) (resp. (3)), label  $ 0  $ the edge  from $ F(10^{n(S)}) $ to $ F(10^{s_{k-2}+1})  $ (resp. $ F(10^{s_{k-1}-g_{l-1}+1})  $) and label  $ 1  $ the edge  from $ F(10^{n(S)}) $ to $ F(1)  $. For a more detailed treatment see \cite{A2}.
 
\subsection{$  \beta$-shifts}
R\'enyi \cite{R} was the first who considered the $ \beta $-shifts. These shifts are symbolic spaces with   rich structures and application in theory and practice. We present here a brief introduction to $ \beta $-shifts from \cite{T2}. For a more detailed treatment, see \cite{B1}.

When $ t  $ is a real number we denote by $ \lfloor t\rfloor $ the largest integer which is smaller than $ t $. Let $ \beta $ a real number greater than $ 1 $. Set 
$$1_{\beta}=a_{1}a_{2}a_{3}\cdots \in \{0,\,1,\ldots,\,\lfloor\beta\rfloor\}^{\mathbb N},$$
where $ a_{1}=\lfloor\beta\rfloor $ and 
$$a_{i}=\lfloor\beta^{i}(1-a_{1}\beta^{-1}-a_{2}\beta^{-2}-\cdots-a_{i-1}\beta^{-i+1})\rfloor$$
for $ i\geq 2 $. The sequence $ 1_{\beta} $ is the expansion of $ 1 $ in the base $ \beta $, that is, $ 1=\sum_{i=1}^{\infty}a_{i}\beta^{-i} $. Let $ \leq $ be the lexiographic ordering of $ (\mathbb N \cup \{0\})^{\mathbb N} $. The sequence $ 1_{\beta} $ has the property that 
\begin{equation}\label{beta1}
\sigma^{k}1_{\beta}\leq 1_{\beta}, \hspace{10mm}k \in \mathbb N,
\end{equation}
where $ \sigma $ denotes the shift on $ (\mathbb N \cup \{0\})^{\mathbb N} $. It is a result of Parry \cite{P}, that this property characterizes the elements of $ (\mathbb N \cup \{0\})^{\mathbb N} $ which are the $ \beta $-expansion of $ 1 $ for some $ \beta>1 $. Furthermore, it follows from \eqref{beta1} that
\begin{equation}\label{beta2}
X_{\beta}=\{x \in \{0,\,1,\ldots,\lfloor\beta\rfloor\}^{\mathbb Z}:\, x_{[i,\,\infty)}\leq 1_{\beta},\,  i \in \mathbb Z\}
\end{equation}
is a shift space of $ \{0,\,1,\ldots,\,\lfloor\beta\rfloor\}^{\mathbb Z} $, called the $ \beta $-shift. 
The $ \beta $-shift is SFT if and only if the $ \beta $-expansion of $ 1 $ is finite and it is sofic if and only if the $ \beta $-expansion of $ 1 $ is eventually periodic \cite{B2}.
Since we are dealing with the case where $\beta\in(1,\,2]$, $a_1=1$ throughout this paper.

\section{Equivalencies between a Beta-shift and its Associate, Sofic Case}

A sliding block code $ \phi:X\rightarrow Y $ is \emph{finite-to-one} if there is
an integer $ M $ such that $ \phi^{-1}(y) $ contains at most $ M $ points for every $ y \in Y $. Shift spaces $ X $ and $ Y $ are \emph{finitely equivalent} if there is an SFT say $ W $ together with finite-to-one factor codes $\phi_{X}:W\rightarrow X $ and $\phi_{Y}:W\rightarrow Y $. We call $ W $ a common extension and $ \phi_{X} $, $ \phi_{Y} $ the \emph{legs}. The triple $ (W, \phi_{X}, \phi_{Y}) $ is a finite equivalence between $ X $ and $ Y $. Call a finite equivalence between sofic shifts in which both legs are right-resolving (resp. right-closing) a \emph{right-resolving finite equivalence} (resp. \emph{right-closing finite equivalence}).

Let $ G $ and $ H $ be two irreducible graphs. Write $ H\preceq G $ if $ X_{H} $ is a right-resolving
factor of $ X_{G} $ and let $ \mathcal R_{G} $ be the collection of graph-isomorphism classes of graphs $ H $
for which $ H\preceq G $. This ordering naturally determines an ordering which
we still call $ \preceq $ on $ \mathcal R_{G} $. Let $M_G$ be the smallest element in the partial ordering $(\mathcal R_{G},\, \preceq)$.

 Now we recall from \cite{LM} how $ M_{G} $ can be constructed. Let $ \mathcal V=\mathcal V(G) $ be the set of vertices of $ G $ and  define a nested
sequence of equivalence relations $ \sim_{n} $ on $ \mathcal V $ for $ n \geq 0 $. The partition of $ \mathcal V $
into $ \sim_{n} $ equivalence classes is denoted by $ \mathcal P_{n} $.
To define $ \sim_{n} $, first let $ I \sim_{0} J $ for all $ I, \,J \in \mathcal V $. For $ n\geq 1 $,
let $ I \sim_{n} J $ if and only if for each class (or atom) $ P \in \mathcal P_{n-1} $ the
total number of edges from $ I $ to vertices in $ P $ equals the total number of edges from
$  J$ to vertices in $  P$.
Note that the partitions  $\mathcal P_{n}$ are nested:
each atom in $  \mathcal P_{n}$ is a union of atoms in $  \mathcal P_{n+1}$.

We have $  \mathcal V$  finite and $  \mathcal P_{n}$   nested; so
 $ \mathcal P_{n} $'s will be equal for all sufficiently large  $  n$. Let $  \mathcal P$ denote this limiting partition. Then  $\mathcal{P}$ will be the set of states of $M_G$. To prevent the confusion between $ M_{G} $ and $ G $, we call a vertex in $ M_{G} $, state and of $ G $ just vertex.

Since for all large enough $  n$, $ \mathcal P=\mathcal P_{n}=\mathcal P_{n+1} $,  for each pair $P, \,Q \in \mathcal P  $ there is  $ k $ such that for each $I \in P  $ there are exactly $  k$ edges in $G$ from $  I$ to vertices in $  Q$. We
then assign $ k$ edges in $M_{G}$ from $  P$ to $  Q$.

Therefore to have $M_G$, for each $n$, we refine the atoms of  $\mathcal{P}_n$ and when $\mathcal{P}_n=\mathcal{P}$, then for each $ P, \,Q \in \mathcal P $ and $ I,\, J \in P $, the total number of paths from $ I $ and $J$  to vertices in $ Q $ and also the length of these paths (with respect to $G$) for both $I$ and $J$ are equal.  Another way to obtain $M_G$  arises from this as follows.

We have $ \mathcal P_{0}=\mathcal V(G) $. Then $ \sim_{1} $ partitions vertices by their out-degrees where for $ X_{\beta} $ and $ X=X(S) $, $ \sim_{1} $ partitions vertices into two atoms, one atom containing the vertices with out-degree one and the other with out-degree two. If $\mathcal{P}\neq \mathcal{P}_1$, for the next step if $ P\in \mathcal{P}_1 $ is refined, then it is turn for $Q$ to be refined where $ Q \in \mathcal{P}_1 $  is any atom having  edges terminating to vertices in $ P $.
\begin{thm}\cite[Theorem 8.4.7]{LM}
Suppose that $ X $ and $ Y $ are irreducible sofic shifts. Let
$ G_{X} $ and $ G_{Y} $ denote the underlying graphs of their Fischer covers respectively. Then $ X $ and $ Y $ are right-resolving finitely equivalent if and
only if $ M_{G_{X}}\cong  M_{G_{Y}}$. Moreover, the common extension can be chosen to be
irreducible.
\end{thm}
A point in $X$ is \emph{doubly transitive} if every word in $\mathcal B(X)$
occurs infinitely often to the left and to the right of its representation. Shift spaces $ X $ and $ Y $ are \emph{almost conjugate} if there is
a shift of finite type $ W $ and 1-1 a.e. factor codes $ \phi_{X}: W \rightarrow X $ and $ \phi_{Y}: W \rightarrow Y $ (1-1 a.e. means that any doubly transitive point has exactly one pre-image). Call an almost conjugacy between sofic shifts in which both legs are right-resolving (resp. right-closing) a \emph{right-resolving almost conjugacy} (resp. \emph{right-closing almost conjugacy}).

Let r-r and r-c stand for right resolving and right closing respectively.
We summarize the relations amongst mentioned properties in the following diagram.
\begin{equation}\label{diagram}
\begin{matrix}
 &&&& \text{conjugacy} \cr
 &&&& \Downarrow \cr
 \text{r-r almost conjugacy} &   \Rightarrow   & \text{r-c almost conjugacy} &   \Rightarrow  & \text{almost conjugacy} \cr
  \Downarrow  &&   \Downarrow   &&   \Downarrow    \cr
 \text{r-r finite equivalence} &   \Rightarrow  & \text{r-c finite equivalence} &   \Rightarrow   & \text{finite equivalence}\cr
\end{matrix}
\end{equation}
There are examples to show that in general
the above properties are different \cite{LM}.
\begin{defn}\label{least period}
Let $w=w_0w_1\ldots w_{p-1}$ be a block of length $ p $. The \emph{least period} of $w$ is the smallest integer $ q $ such that $w=(w_0w_1\ldots w_{q-1})^m$ where $ m=\frac{p}{q} $ must be an integer. The block $ w $ is \emph{primitive} if its least period equals its length $ p $.
\end{defn}
Now we set up to picture out the graph $ M_{G} $ of $ X(S) $.
 First suppose $|S|<\infty $.
Let $ S=\{s_{0},s_{1},\ldots,s_{k-1}\}\subseteq\mathbb N_{0} $, $k>1$ and
\begin{equation}\label{D(S)}
  \mathcal D(S)=d_{1}d_{2}\cdots d_{k-2}(d_{k-1}+s_{0}+1)
\end{equation}
where $ d_{i}=s_{i}-s_{i-1} $, $ 1\leq i \leq k-1 $. Note that if $ I, J \in \mathcal V(G) $ are in the same state of $ M(G) $, then both of $ I $ and $ J $ have the same out-degree which is one or two. Also the out-degree of any vertex $ F(10^{s_{i}}) $, $ 0\leq i\leq k-1 $ is two except that the last one. Hence $ d_{i} $, $ 1\leq i\leq k-2 $ measures the distance between any two vertices with out-degree two. 

To pick the next vertex after $ F(10^{s_{k-2}}) $ with out-degree two we continue to right to $ F(10^{s_{k-1}}) $ and then along the graph to $ F(1) $ and then again to right to $ F(10^{s_{0}}) $ that is after $ d_{k-1}+s_{0}+1 $ steps.
\begin{thm}\label{finite}
Let $ |S|<\infty $. Then $\mathcal D(S)  $ is primitive if and only if $ M_{G}\cong G $.
\end{thm}
\begin{proof}
Suppose $\mathcal D(S)  $ is not primitive. Let $ \mathcal V=\mathcal V\mathcal(M_{G}) $ be the set of states of $ M_{G} $. Then by the Fischer cover of $ X(S) $, each state in $ M_{G} $ consists of $ m=\frac{|S|-1}{q-1} $ vertices of graph $ G $ where $ q-1 $ is the least period $\mathcal D(S)  $ and $ |\mathcal V|=\sum_{i=1}^{q-1}d_{i}=s_{q-1}-s_{0}$. In fact, if $\mathcal V=\{P_{i}: 0\leq i \leq s_{q-1}-s_{0}-1\}  $, then
$$P_{i}=\{F(10^{s_{0}+i}),\,F(10^{s_{0}+i+|\mathcal V|}),\ldots,F(10^{s_{0}+i+(m-1)|\mathcal V| \mod u})\}$$
where $ u=s_{k-1}+1 $. Since $ |\mathcal V|=s_{q-1}-s_{0}<s_{k-1}+1=|\mathcal V(G)|  $,  $ M_{G}\not\cong G $.

Now suppose $ M_{G}\not\cong G $. So there are at least two different vertices of $G$ say $ I=F(10^{p}) $ and $ J=F(10^{q}) $  such that $ I$ and $ J $ are in the same state of $ M_{G} $. Assume $p<q$. There exists an edge from $ I $ (resp. $J$) to $ F(10^{(p+1)})$ (resp. $F(10^{(q+1)\mod u})$). Therefore, by the fact that  $ I $ and $ J $ are equivalent, we have that the vertices $ F(10^{(p+1)})$ and  $F(10^{(q+1)\mod u}) $ are equivalent. By the same reasoning, for each $ i\geq 2$, $ F(10^{(p+i)\mod u})$ and  $F(10^{(q+i)\mod u}) $ are equivalent. So $\mathcal D(S)  $ is not primitive.
\end{proof}
\begin{thm}\label{infinite}
Let $ X(S) $ be a sofic shift with $|S|=\infty$ and the Fischer cover $ \mathcal G=(G,\mathcal L) $. Then $ M_{G}\cong G $.
\end{thm}
\begin{proof}
We consider the three cases appearing for $ |S|=\infty $ in subsection \ref{subs:S-gap}. We claim that the last vertex $ F(10^{n(S)}) $ is not equivalent with any other vertex. Otherwise, we will show that at least one of $ k $ or $ l $ will not be the least integer in \eqref{Delta}. So the state of $ M_{G} $  containing this last vertex, contains only this vertex which this in turn implies that other states of $ M_{G} $ also have one vertex. So $ M_{G}\cong G $.

 We prove our claim for the most involved case, that is case (3). First suppose there is a vertex 
\begin{equation}\label{A}
v_{0}=F(10^{t_{0}})\sim F(10^{n(S)}), \hspace{10mm}s_{k-1}-g_{l-1}+1 \leq t_{0}< n(S). 
\end{equation}
Without loss of generality assume that this $ t_{0} $ is the largest integer with this property. Recall that there is an edge from
 $ F(10^{n(S)}) $ to $ F(10^{s_{k-1}-g_{l-1}+1}) $; so it is convenient to set $ t_{1}:=n(S) $, $ t_{1}+1:=s_{k-1}-g_{l-1}+1 $ and $ v_{1}:=F(10^{n(S)}) $. By \eqref{A}, $ v_{2}:=F(10^{t_{1}+1})\sim F(10^{t_{0}+1}) $ and moving horizontally to right $
 F(10^{t_{1}+i})\sim F(10^{t_{0}+i}) $, $ 2\leq i\leq t_{1}-t_{0}$.
 Moreover, none of $ F(10^{t_{0}+i})$ will be equivalent to $ v_{0} $, for this will violate the way we have picked $ t_{0} $. If $ v_{2}\sim v_{0} $ we are done, for then $ l $ will not be the least integer.
Observe that  there are only finitely many vertices; therefore,  there must be $v_{i}\not\sim v_{0} $, $ 2\leq i < p $ and $ v_{p}\sim v_{0} $. 
Applying the same reasoning,
we deduce that again $ l $ is not the least integer.

If $ F(10^{n(S)})  $ is not equivalent to any vertex $F(10^{t}) $ for $ s_{k-1}-g_{l-1}+1 \leq t < n(S) $, it will be equivalent to $ F(10^{s_{k-1}-g_{l-1}}) $. This implies $ k $ is not the least integer.
\end{proof}
Theorems \ref{finite} and \ref{infinite} imply the following.
\begin{cor}\label{MG}
Let $ X(S) $ be a sofic shift with the Fischer cover $ \mathcal G=(G,\mathcal L) $. Then any state of $ M_{G} $ has the same number of vertices of $ G $.
\end{cor}
When $ |S|<\infty $, there may be cases with $ M_{G}\not\cong G $. The difference with $ |S|=\infty $ is that for $ |S|<\infty $, the last vertex $F(10^{n(S)})$ has always out-degree one with label $ 1 $ while for $ |S|=\infty $, the label of edge starting from the vertex with out-degree one is $ 0 $.

Now let $ X $ be a sofic shift with the Fischer cover $ \mathcal G=(G, \mathcal L) $. Then by definition, $ \mathcal L_{\infty} $ is right-resolving and also it is almost invertible \cite[Proposition 9.1.6]{LM}. So
\begin{lem}\label{almost invertible} 
 If $ X $ and $ Y $ are sofic with Fischer covers $ \mathcal G_{X}=(G_{X}, \mathcal L_{X}) $ and $ \mathcal G_{Y}=(G_{Y}, \mathcal L_{Y}) $ respectively, such that $G_{X} \cong G_{Y}  $, then $ X $ and $ Y $ will be right-resolving almost conjugate with legs $ {\mathcal L_{X}}_{\infty} $ and $ {\mathcal L_{Y}}_{\infty} $.
\end{lem}
\begin{thm}\label{finite equivalent}
Let $ X_{\beta}$ be a sofic $ \beta $-shift for $ \beta \in (1,\,2] $. Then there is $ S \subseteq \mathbb N_{0} $, determined in terms of
 coefficients of $1_\beta$, such that $ X_{\beta} $ and $ X(S) $ are right-resolving almost conjugate. Moreover, if $ X_{\beta} $ is SFT, then  $ X(S) $ can be chosen to be  conjugate to $ X_{\beta} $.
\end{thm}
\begin{proof}
For a given sofic $ \beta $-shift, $ \beta \in (1,\,2] $, we claim that there is 
 $S\subseteq\mathbb N_{0}$  such that the $S$-gap shift $ X(S) $ and $ X_{\beta} $ have the same underlying graph for their Fischer covers. Then by Lemma \ref{almost invertible}, $ X_{\beta} $ and $ X(S) $ will be right-resolving almost conjugate. 

Let $ 1_{\beta}=a_{1}a_{2}\cdots a_{n}(a_{n+1}\cdots a_{n+p})^{\infty} $ and
$\{i_{1},\,i_{2},\ldots,\,i_{t}\}\subseteq\{1,\,2,\ldots,\,n\}$ where $a_{i_{v}}=1$ for $1\leq v\leq t  $. 
Note  that $i_1$ is always 1.
Similarly, let $\{j_{1},\,j_{2},\ldots,\,j_{u}\}\subseteq\{n+1,\,\ldots,\,n+p\}$ where $a_{j_{w}}=1$ for $1\leq w \leq u$.
We consider two cases:
\begin{enumerate}
\item
  $ X_{\beta}$ is SFT. In this case $(a_{n+1}\cdots a_{n+p})^{\infty}=0^{\infty}  $ and $ a_{n}=1 $. So $ i_{t}=n $ and $X(S)$ with
  \begin{equation}\label{Sgap-unique}
S=\{0,\,i_{2}-1,\,\ldots,\,i_{t-1}-1,\,i_{t}-1\}.
\end{equation}
is the required $S$-gap shift as has been claimed. Since both  $ X_{\beta} $ and $ X(S) $ are SFT with the same underlying graph $G$ for their  Fischer covers, they are both conjugate to $ X_{G} $ \cite[Theorem 3.4.17]{LM},  and so conjugate to each other.
 
\item
 $ X_{\beta} $ is strictly sofic. Then $(a_{n+1}\cdots a_{n+p})^{\infty}\neq 0^{\infty}  $. Relabel any edge on $ G_{\beta} $ ending at the first vertex for $ 1 $ and other edges for $ 0 $. Shift space corresponding to this labeling is an $ S $-gap shift where
\begin{equation}\label{S}
S=\{0,\,i_{2}-1,\,\ldots,\,i_{t}-1,\,j_{1}-1,\,\ldots,\,j_{u}-1,\,j_{1}+p-1,\ldots\}.
\end{equation}
(Observe that then
\begin{equation}\label{Delta 2}
\Delta(S)=\{0,\,i_{2}-1,\ldots,\,i_{t}-i_{t-1},\,j_{1}-i_{t},\,\overline{j_{2}-j_{1},\ldots,\,j_{u}-j_{u-1},\,j_{1}-j_{u}+p}\}
\end{equation}
Which shows that $ X(S) $ is sofic \cite[Theorem 3.4]{A1}.)

Rewrite $ \Delta(S) $ in \eqref{Delta 2} as
$$\Delta(S)=\{0,d_{1},\ldots,d_{t},\overline{g_{0},\ldots,g_{u-1}}\}.$$
We claim that $ \mathcal G_{S} =(G_{S},\,\mathcal L_{S})$ is follower-separated. Otherwise, there are two cases.
\begin{enumerate}
\item
There is $ 1\leq i \leq t $ such that $ d_{t+1-j}=g_{u-j} $, $ 1\leq j \leq i $. Then $ \mathcal G_{\beta}=(G_{\beta}, \mathcal L_{\beta}) $ is not follower-separated and so it is not the Fischer cover of $ X_{\beta}$ which is absurd.
\item
$ g_{0}g_{2}\cdots g_{u-1} $ is not primitive. This implies that $ a_{n+1}\cdots a_{n+p}$ is not primitive which is again absurd.
\end{enumerate}
This establishes the claim and  $S$ is completely determined.
\end{enumerate}
\end{proof}
Now the following is immediate.
\begin{cor}\label{MG beta}
Let $ X_{\beta} $ be a sofic shift whose underlying graph of its Fischer cover is $G$. Then $ M_{G}\cong G $.
\end{cor}
\begin{proof}
Suppose $ M_{G}\not\cong G $. For this $ X_{\beta} $, find the $S$-gap shift satisfying the conclusion of Theorem \ref{finite equivalent}. Then by Theorem \ref{infinite}, this $X(S)$ (as well as our $X_\beta$) must be SFT and $ \mathcal D(S) $ is not primitive. But this will not allow to have \eqref{beta1} which is a necessary condition.
\end{proof}
\begin{lem}\label{case one}
Let $ |S|=\infty $ and $ X(S) $ be a sofic shift satisfying (1a) in Subsection \ref{subs:S-gap}. Then there does not exist any $ \beta $-shift being right-resolving finite equivalent with $ X(S) $.
\end{lem}
\begin{proof}
Suppose there is some $ \beta \in (1,\,2] $ such that $ X(S) $ and $ X_{\beta} $ are right-resolving finite equivalent and
 $ \mathcal G_{S}=(G_{S},\mathcal L_{S}) $ and $ \mathcal G_{\beta}=(G_{\beta},\mathcal L_{\beta}) $ are the Fischer covers of $ X(S) $ and $ X_{\beta}$ respectively.
 By Theorem \ref{infinite} and Corollary \ref{MG beta}, $ G_{S}\cong G_{\beta} $. Then $ G_{\beta} $ is the underlying graph of 
$ \mathcal G_{S} $ and $ 1_{\beta}=(a_{1}a_{2}\cdots a_{n})^{\infty} $.

Now by hypothesis, $ g_{l-1}=1 $, so $ 1\notin S $ and this implies that $ a_{2}=0 $ while $ a_{1}=a_{n}=1 $. This means $ (a_{1}a_{2}\cdots a_{n})^{\infty} $ does not satisfy \eqref{beta1} and we are done.
\end{proof}
Let $ X(S) $ be an $ S $-gap shift where $ s_{0}=0 $ and $ d_{i}=s_{i}-s_{i-1} $, $ i \in \mathbb N $ and also $\mathcal{D}(S)$  as \eqref{D(S)}. Define
\begin{equation}\label{inverse Parry}
 d_{1}d_{2}d_{3}\cdots =\left\{
\begin{tabular}{cll}
$(d_{1}d_{2}\cdots(d_{k-1}+1))^{\mathbb{N}}=(\mathcal D(S))^{\mathbb N},$&$|S|=k;$ \\
$d_{1}d_{2}\cdots,$& $|S|=\infty.$ \\
\end{tabular}
\right.
\end{equation}
\begin{thm}\label{inverse}
Suppose $ X(S) $ is a sofic shift.  Then $ X(S) $ is right-resolving almost conjugate to a $ \beta $-shift if and only if
\begin{equation}\label{varepsilon}
d_{n}d_{n+1}\cdots\geq d_{1}d_{2}\cdots
\end{equation}
for all $ n\geq 1 $. 
\end{thm}
\begin{proof}
 Let $ \beta \in (1,\,2] $ with $ 1_{\beta}=a_{1}a_{2}\cdots $ be so that $ X(S) $ and 
$ X_{\beta} $ are right-resolving almost conjugate. This means they are
 right-resolving finite equivalent. First suppose $ M_{G_{S}}\cong G_{S} $. By Corollary \ref{MG beta}, $G_{S}\cong G_{\beta} $ and so
 \eqref{varepsilon} follows from the fact that $ a_{1}a_{2}\cdots $ satisfies \eqref{beta1}.

If $ M_{G_{S}}\not\cong G_{S} $, then by Theorems \ref{finite} and \ref{infinite}, $ |S|<\infty $. So $ X_{\beta} $ is right-resolving finite
 equivalent to $ X(S') $ with $ S'=\{0,s_{1},\ldots,(s_{q-1}-1)\} $ and $ \mathcal D(S)=\mathcal D(S')^{m} $ where $ m=\frac{|S|-1}{q-1} $ as in the proof of Theorem \ref{finite}. Moreover, $M_{G_{S'}}\cong G_{S'}  $ which gives again $d'_{n}d'_{n+1}\cdots\geq d'_{1}d'_{2}\cdots$ for all $ n\geq 1 $. Now this fact reflects to $ \mathcal D(S) $ and \eqref{varepsilon} holds.

To prove the sufficiency suppose $ \mathcal G_{S}=(G_{S},\mathcal L_{S}) $ is the Fischer cover of $ X(S) $ and $ \mathcal V=\mathcal V(G_{S}) $  the set of vertices of $ G $. Relabel $ G_{S} $ by labeling $ 0 $ any edge terminating at vertex $ F(1) $ and any edge whose initial vertex has out-degree $ 1 $, and assign $ 1 $  all other edges. 
 
 Recall that we have lined up the vertices horizontally from $F(1)$ in left to $F(10^{n(S)})$ on right. First let $ |S|<\infty $ and $ a_{1}a_{2}\cdots a_{n(S)} $ be the assigned label of the horizontal path from $ F(1) $ to the last vertex with $ a_{i}=0 $ or $ 1 $ as determined above. Then \eqref{varepsilon} implies that $ a_{1}a_{2}\cdots a_{n(S)}1 $ is the $ \beta $-expansion of $ 1 $ for
 some $ \beta \in (1,\,2] $ and $ \mathcal G_{\beta}$ is the Fischer cover of $ X_{\beta} $. 
 
 When $ |S|=\infty $, assign the label $ a_{1}a_{2}\cdots a_{n(S)} $ to the horizontal path from $ F(1) $ to the last vertex and label $ a_{n(S)+1} $ to the edge starting from $ F(10^{n(S)}) $ and terminating at $ F(10^{n(S)+1}) $. Again \eqref{varepsilon} implies that $ a_{1}a_{2}\cdots a_{n}(a_{n+1}\cdots a_{n(S)+1})^{\infty} $ is the $ \beta $-expansion of $ 1 $ for
 some $ \beta \in (1,\,2] $ where the index $ n $ depends on $ S $. Then $ \mathcal G_{\beta}$ is the Fischer cover of $ X_{\beta} $ (one needs similar arguments as in the proof of
 Theorem \ref{finite equivalent} to see this fact). So Lemma \ref{almost invertible} implies that $ X(S) $ and $ X_{\beta} $ are right-resolving almost conjugate.
\end{proof}
Let $ S=\{n,\,n+1,\ldots\} $, $ n \geq 2 $ and $ 1_{\beta}=10^{n}1 $. Then $ X_{S} $ and $ X_{\beta} $ are right-resolving almost conjugate;
 however, the condition \eqref{varepsilon} does not hold. This is not a contradiction, for we are  considering $ X(S) $ up to conjugacy and in this exceptional case, we consider $X(S')=$ASS$(X_\beta)$ where $S'=\{0,n\}$.
\begin{rem}
  $ X_{\beta} $ can be explicitly determined  in terms of $ S $. If $ S=\{0,\,s_{1},\ldots,\,s_{k-1}\} $, then it is sufficient to set $ 1_{\beta}=a_{1}a_{2}\cdots a_{s_{k-1}+1} $ such that 
$ a_{1}=a_{s_{i}+1}=1 $, $ 1\leq i \leq k-1$. When $ |S|=\infty $, different cases of Subsection \ref{subs:S-gap} must be considered. Case (1a) has been ruled out by Lemma \ref{case one}, so other cases will be considered. 
\begin{list}{•}{•}
\item{(1b)}
if $ k=1 $ and $ g_{l-1}>1 $, then $ F(10^{g})=F(1) $. So $ 1_{\beta}=a_{1}a_{2}\cdots a_{g}$ such that $ a_{s_{i}+1}=1 $, $ 0\leq i \leq l-1$.
\item{(2)}
if $ k\neq1 $ and $ g_{l-1}>d_{k-1} $, then $ F(10^{g+s_{k-2}+1})=F(10^{s_{k-2}+1}) $. So $ 1_{\beta}=a_{1}a_{2}\cdots a_{s_{k-2}+1}(a_{s_{k-2}+2}\cdots a_{g+s_{k-2}+1})^{\infty}$ which $ a_{s_{i}+1}=1 $, $ 0\leq i \leq k+l-2$.
\item{(3)}
if $ g_{l-1}\leq d_{k-1} $, then $ F(10^{s_{k+l-2}+1})=F(10^{s_{k-1}-g_{l-1}+1}) $. So $$ 1_{\beta}=a_{1}a_{2}\cdots a_{s_{k-1}-g_{l-1}+1}(a_{s_{k-1}-g_{l-1}+2}\cdots a_{s_{k+l-2}+1})^{\infty}$$
 which $ a_{s_{i}+1}=1 $, $ 0\leq i \leq k+l-2$ and $ a_{s_{k+l-2}+1}=1 $.
\end{list}
\end{rem}
Now we show that the conclusion of Theorem \ref{finite equivalent} about conjugacy is not true in non-SFT case. Recall that when $  X$ is a shift space with non-wandering part $ R(X) $, we can
consider the shift space
$$\partial X =\{x \in R(X) : x \mbox{ contains no words that are synchronizing for } R(X)\};$$
which is called the \emph{derived shift space} of $  X$. The derived shift space  is a conjugacy invariant.
\begin{thm}\label{non-conjugacy}
A non-SFT  $ \beta $-shift is not conjugate to an $S$-gap shift for any $S\subseteq \mathbb N_{0}$. 
\end{thm}
\begin{proof}
All the $S$-gap shifts are synchronized; therefore, a possible conjugacy happens between  synchronized $\beta$ and  $S$-gap shifts and so we assume that our non-SFT $\beta$-shift is  synchronized.

Suppose that there is $ S\subseteq \mathbb N_{0}$ such that $ \varphi: X(S) \rightarrow X_{\beta}$ is a conjugacy map. Then
by \cite[Proposition 4.5]{T3},
 we must have $ \varphi(\partial X(S))=\partial X_\beta $.  Since $ 1 $ is a synchronizing word for any  $ S $-gap shift, and $X(S)$ is not SFT, $ \partial X(S)=\{0^{\infty}\} $ (for an SFT $S$-gap shift, $ \partial X(S)=\emptyset $).
To prove the theorem, we show that 
\begin{equation}\label{namosavi}
\varphi( \{0^{\infty}\})\neq \partial X_\beta.
\end{equation}
Recall that the $ \omega $-limit set of the sequence $ 1_{\beta} $ under the shift map  is the derived shift space $ \partial X_{\beta} $ of $ X_{\beta} $ \cite[Theorem 2.8]{T2}.
First assume that $ X_{\beta} $ has specification property. Then there exists some $ n\geq 0 $ such that $ 0^{n} $ is not a factor of $ 1_{\beta} $ \cite{B3}. So $ 0^{n} $ is a synchronizing word for $ X_{\beta} $ \cite[Proposition 2.5.2]{B3} and $ 0^{\infty}\not\in \partial X_{\beta} $. Therefore, $ \partial X_{\beta} \cap P_{1}(X_{\beta})=\emptyset  $  while $ \varphi(0^{\infty})\in P_{1}(X_{\beta})$ and $ \varphi(0^{\infty})\in \varphi(\partial X(S))=\partial X_{\beta}$ and \eqref{namosavi} holds.

If $ X_{\beta} $ does not have specification, then $ \{0^{\infty}, 10^{\infty}\}\subseteq \omega(1_{\beta})=\partial X_{\beta} $ and again \eqref{namosavi} holds.
\end{proof}
\begin{cor}\label{unique}
Let $ X_{\beta} $ be SFT and $X(S_0)$ the unique $S$-gap shift conjugate to $X_\beta$ (Theorem \ref{finite equivalent}).
Then $ X_{\beta} $ is   
\begin{enumerate}
\item
right-resolving
almost conjugate to $X(S_0)$,
\item
right-resolving
finite equivalent 
 to infinitely many $ S $-gap shifts $ (X(S_{n}))_{n \in \mathbb N} $
 with $ \mathcal D(S_{n})=(\mathcal D(S_{0}))^{n+1} $, $ n \in \mathbb{N} $,
\item
right-resolving
almost conjugate to a unique strictly sofic $S$-gap shift.
\end{enumerate}
 If $ X_{\beta} $ is strictly sofic, then it is
right-resolving almost conjugate to a unique $ S $-gap shift. 
\end{cor}
\begin{proof}
Let $X_\beta$ be SFT and let $ 1_{\beta}=a_{1}a_{2}\cdots a_{n-1}a_{n} $ and $$\{i_{1},\,i_{2},\ldots,\,i_{t}\}\subseteq\{1,\,2,\ldots,\,n\}$$ where $ a_{i_{j}}=1 $, $ 1\leq j \leq t $. We will relable the  Fischer cover of $X_\beta$ for possible presentation of an $S$-gap shift. 

One  of such SFT $ S $-gap shifts is $X(S_0)$ characterized in the proof of Theorem \ref{finite equivalent}. By that theorem, $ X_{\beta} $
 and $ X(S_{0}) $ are right-resolving almost conjugate and conjugate which gives (1). For (2) relabel $ \Delta(S_{0})=\{0,\,s_{2}-1,\,s_{3}-s_{2},\ldots,\,s_{t}-s_{t-1}\} $ as $  \Delta(S_{0})=\{0,\,d_{1},\ldots,\,d_{t-1}\} $ and observe that $ \mathcal D(S_{0})=d_{1}\cdots d_{t-2}(d_{t-1}+1) $.
 Set $ \mathcal H(S_0):=\{0,\,d_{1},\,d_{1}+d_{2},\ldots,\,\sum_{i=1}^{t-1}d_{i}\} $ and let
$$S_{1}=\left(S_{0}\setminus \{i_{t}-1\}\right)\cup (i_{t}+\mathcal H(S_{0})).$$
Then $ \mathcal D(S_{1})=(\mathcal D(S_{0}))^{2} $  is not primitive and we have $ M_{G_{S_{1}}}\cong M_{G_{S_{0}}}$. 

Now for $j\in\mathbb{N}$, let $ s_{i_{j}}=\max \{s:\, s \in S_{j-1}\} $ and use an induction argument to see that for 
\begin{equation}\label{S_j}
S_{j}=(S_{j-1}\setminus \{s_{i_{j}}\})\cup ((s_{i_{j}}+1)+\mathcal H(S_{0})),
\end{equation}
$ \mathcal D(S_{j})=(\mathcal D(S_{0}))^{j+1} $ and $ M_{G_{S_{j}}}\cong M_{G_{S_{0}}}$.

To prove (3) note that there is an strictly sofic $S$-gap shift with $ k=1 $ and $ g_{l-1}>1 $ as in Subsection \ref{Delta} where 
$ S=\{0,\,i_{2}-1,\ldots,\,i_{t-1}-1,\,i_{t},\,i_{t}+i_{2}-1,\ldots\} $. The element $i_{t}  $ appears in $ S $ because the edge starting from the last vertex and terminating at the first vertex is labeled $ 0 $. In fact
$$\Delta(S)=\{0,\overline{i_{2}-1,i_{3}-i_{2},\ldots,i_{t-1}-i_{t-2},i_{t}-i_{t-1}+1}\}.$$
Hence $ X_{\beta} $ and $ X(S) $ have the same underlying graph for their Fischer covers and by  Lemma \ref{almost invertible}, they are right-resolving almost conjugate.

If there is another $ S $-gap shift such that $ X_{\beta} $ and $ X(S) $ are right-resolving finite equivalent, then $ M_{G_{\beta}}\cong M_{G_{S}}$ and so $ M_{G_{S_{0}}}\cong M_{G_{S}}$. 
Now Theorems \ref{finite} and \ref{infinite} imply that $ |S|<\infty $ and $ \mathcal D(S) $ is not primitive which in turn implies that $ \mathcal D(S)=(\mathcal D(S_{0}))^{m} $ for some $ m \in \mathbb N $. Therefore, $ S=S_{m-1} $ as defined in \eqref{S_j}. 

Now suppose $ X_{\beta} $ is strictly sofic. A typical Fischer cover of $ X_{\beta} $ has been shown in Figure 3.1. The existence of loop in the first vertex from left implies that it is the vertex $ F(1) $ in the Fischer cover of $ S $-gap shift. By Fischer cover of $ S $-gap shifts \cite{A2}, there is only one $ X(S) $ with the Fischer cover as appears in Figure 3.1. 
\end{proof}
\begin{center}\label{graph} 
\begin{tikzpicture}[->,>=stealth',shorten >=1pt,auto,node distance=1.5cm,
                    thick,main node/.style={rectangle,draw,font=\sffamily\small\bfseries}]

  \node[main node] (1) {};
\node[] (2) [ right of=1] {$\ldots$};
\node[main node] (3) [right of=2] {};
\node[] (4) [ right of=3] {$\ldots$};
\node[main node] (5) [right of=4] {};
  \node[] (6) [ right of=5] {$\ldots$};
  \node[main node] (7) [right of=6] {};

 \path[every node/.style={font=\sffamily\small}]
    (1) edge node [right] { $\hspace{-6mm}\ ^{\ ^{a_{1}}}$} (2) 
    edge [loop above] node {$\hspace{-2mm}\ _{a_{1}-1}$} (1)
    (2) edge node [right] {$\hspace{-6mm}\ ^{\ ^{a_{i-1}}}$} (3) 
    (3) edge node [right] {$\hspace{-6mm}\ ^{\ ^{a_{i}}}$} (4)
     edge [bend left] node[right] {$\hspace{-6mm}\ ^{\ ^{a_{i}-1}}$} (1)
(4) edge node [right] {$\hspace{-6mm}\ ^{\ ^{a_{n}}}$} (5)
(5) edge node [right] {$\hspace{-6.5mm}\ ^{\ ^{a_{n+1}}}$} (6)
edge [bend left] node[right] {$\hspace{-6mm} \ ^{\ ^{a_{n+1}-1}} $} (1)
(6) edge node [right] {$ \hspace{-9.5mm}\ ^{\ ^{a_{n+p-1}}} $} (7)
(7) edge [bend right] node[right] {$\hspace{-6mm} \ ^{\ ^{a_{n+p}}} $} (5) 
 edge [bend left] node[right] {$\hspace{-6mm}\ ^{\ ^{a_{n+p}-1}}$} (1);
\end{tikzpicture}
\end{center}
\begin{center}
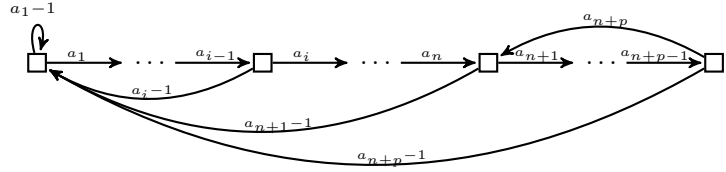
Figure 3.1: A typical Fischer cover of an strictly sofic $ \beta $-shift for $ 1_{\beta}=a_{1}a_{2}\cdots a_{n}(a_{n+1}\cdots a_{n+p})^{\infty} $, $ \beta \in (1,\,2] $. The edges 
heading to $\alpha_1$ exist if $a_i=1$.
\end{center}
It is not hard to see that for a non-sofic $\beta$ shift, like strictly sofic case, there exists a unique $S$-gap satisfying \eqref{varepsilon} and having the same underlying graph for its Fischer cover as our $X_\beta$.
This motivates the following definition.
\begin{defn}\label{defn:associate}
We say that $X(S)$ is the associated $S$-gap shift to a $\beta$-shift and is denoted by ASS$(X_\beta)$, $\beta\in[1,\,2)$
if
\begin{itemize}
\item
$X_\beta$ is SFT and conjugate to $X(S)$, or
\item
$X_\beta$ is not SFT and has the same
underlying graph for its Fischer cover as  $X(S)$.
\end{itemize}
Similarly, for $X(S)$ satisfying \eqref{varepsilon}, a unique $X_\beta$ exists such that $X(S)=$ASS$(X_\beta)$. 
This $X_\beta$ is called the associated $\beta$-shift to $X(S)$ and is denoted by ASS$(X(S))$.
\end{defn}
 Therefore, ASS(ASS$(X_\beta))=X_\beta$ and ASS(ASS$(X(S)))=X(S)$ for $\beta\in[1,\,2)$ and $S$ satisfying \eqref{varepsilon}.
\begin{rem}
 $X_{\beta} $ and ASS$( X_\beta) $ in Theorem \ref{finite equivalent} have all
 equivalencies given in
 diagram \eqref{diagram} when they are both SFT and all
except conjugacy when they are strictly sofic.
\end{rem}

\section{Common Properties between a Beta-shift and its Associate, Non-sofic Case}
By \cite[Theorem 4.22]{J}, for every $ \beta>1 $
 there exists $ 1<\beta'<2 $ such that $ X_{\beta} $ and $ X(\beta') $ are flow equivalent.
But any two flow equivalent shift spaces have the same Bowen-Franks groups. Therefore, by  Theorem \ref{finite equivalent} and \cite[Theorems 3.1 and 3.2]{A2} which gives a complete account of the Bowen-Franks group of $S$-gap shifts, we have also a complete characterization of such groups for $\beta$-gap shifts for  $ \beta>1 $.

Also when $ X_{\beta} $ is 
 sofic shift with  Fischer cover  $ \mathcal G_{\beta}=(G_{\beta},\mathcal L_{\beta}) $, by Theorems \ref{finite equivalent} and \cite[Theorem 2.2]{A2} we can determine the characteristic polynomial of $ G_{\beta} $ and so we have all eigenvalues.

Now we consider non-sofic $ \beta $ and $ S $-gap shifts. By Theorem \ref{non-conjugacy}, no conjugacy occurs between a $\beta$-shift and any $S$-gap shift and hence we set up to look for other possible equivalencies, and in particular, similar to those in diagram \eqref{diagram}.

First recall
 that any $ \beta $-shift is \emph{half-synchronized} \cite{FF} whereas any $ S $-gap shift is synchronized. So the Fischer covers $ \mathcal G_{\beta}=(G_{\beta},\mathcal L_{\beta}) $ and $ \mathcal G_{S}=(G_{S},\mathcal L_{S}) $ exist (see Figure 4.1).

Let $1_{\beta}=a_{1}a_{2}\cdots$. Relabel $ G_{\beta} $ by labeling $ 1 $ any edge terminating at vertex $ \alpha_{1} $ and $ 0 $ all other edges to get an $ S $-gap shift with the same underlying graph as $ X_{\beta} $. Also using Theorem \ref{inverse},
one can relabel a non-sofic $ S $-gap shift satisfying \eqref{varepsilon}  to obtain a $ \beta $-shift. In both cases, $ X_{S}=ASS(X_{\beta}) $. 
\begin{center}\label{graph1} 
\begin{tikzpicture}[->,>=stealth',shorten >=1pt,auto,node distance=1.5cm,
                    thick,main node/.style={rectangle,draw,font=\sffamily\small\bfseries}]

  \node[main node](1) {$ \alpha_{1} $};
\node[] (2) [ right of=1] {$\ldots$};
\node[main node] (3) [right of=2] {$ \alpha_{i} $};
\node[] (4) [ right of=3] {$\ldots$};
\node[main node] (5) [right of=4] {$ \alpha_{n} $};
  \node[] (6) [ right of=5] {$\ldots$};

 \path[every node/.style={font=\sffamily\small}]
    (1) edge node [right] { $\hspace{-6mm}\ ^{\ ^{a_{1}}}$} (2) 
    edge [loop above] node {$\hspace{-2mm}\ ^{a_{1}-1}$} (1)
    (2) edge node [right] {$\hspace{-7mm}\ ^{\ ^{a_{i-1}}}$} (3) 
    (3) edge node [right] {$\hspace{-6mm}\ ^{\ ^{a_{i}}}$} (4)
     edge [bend left] node[right] {$\hspace{-6mm}\ ^{\ ^{a_{i}-1}}$} (1)
(4) edge node [right] {$\hspace{-7mm}\ ^{\ ^{a_{n}-1}}$} (5)
(5) edge node [right] {$\hspace{-7mm}\ ^{\ ^{a_{n}}}$} (6)
edge [bend left] node[right] {$\hspace{-6mm} \ ^{\ ^{a_{n}-1}} $} (1);


\end{tikzpicture}
\end{center}
\begin{center}
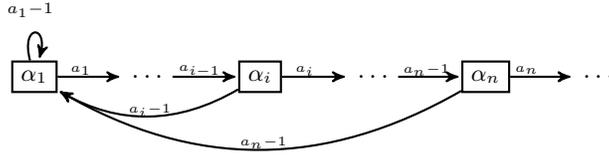
Figure 4.1: A typical Fischer cover of a non-sofic $ \beta $-shift for $ 1_{\beta}=a_{1}a_{2}\cdots $, $ \beta \in (1,\,2] $. The edges 
heading to $\alpha_1$ exist if $a_i=1$.
This cover can be relabeled to give the Fischer cover of  ASS$(X_\beta)$.
\end{center}
\begin{thm}\label{entropy} 
$ h(X_{\beta})=h($ASS$(X_\beta)), \beta\in (1,\,2]$. 
\end{thm}
\begin{proof}
Entropy is an invariant for all the properties given in diagram \eqref{diagram}. So when $ X_{\beta} $ is sofic, the proof is obvious (Theorem \ref{finite equivalent}). 

Now let $ X_{\beta} $ be a non-sofic shift and  let $
 1_{\beta}=a_{1}a_{2}\cdots $. We have $ a_{i}=1 $ if and only if $ i-1 \in S $. But for  $ 1_{\beta}=\sum_{i=1}^{\infty}a_{i}\beta^{-i} $, $ h(X_{\beta})=\log \beta $ and $ h(X(S))=\log \lambda $ where
$\lambda$ is a nonnegative solution of  $\sum_{n\in S}x^{-(n+1)}=1$ \cite{S}. So $ h(X_{\beta})=h($ASS$(X_\beta))$.
\end{proof}
Let $  X$ and $Y  $ be two coded systems. Then there is a coded system $  Z$ which
factors onto $  X$ and $Y  $ with entropy-preserving maps if and only if $ h(X)=h(Y) $. In particular, $  Z$ can be chosen to be an almost Markov synchronized system \cite[Theorem 2.1]{F}. So Theorem \ref{entropy} implies the following.
\begin{cor}\label{almost Markov} 
$ X_{\beta} $ and ASS$(X_{\beta}) $ have a common extension 
which is an
almost Markov synchronized system whose maps are entropy-preserving.
\end{cor}
For a dynamical system $(X,T)$, let $p_{n}$ be the number of periodic points in $X$ having period $n$. When $p_{n} < \infty$, the zeta function $\zeta_{T}(t)$ is defined as
$$
\zeta_{T}(t)=\exp\left(\sum_{n=1}^{\infty}\frac{p_{n}}{n}t^{n}\right).
$$
The zeta functions of $ \beta $-shifts have been determined in \cite{FLP}. 
Here we will give the zeta function of $\zeta_{\sigma_S}$ in terms of $\zeta_{\sigma_\beta}$ where $X(S)=$ASS$(X_\beta)$.

Let $ X_{\beta} $ be sofic and  $1_{\beta}=a_{1}a_{2}\cdots a_{n}(a_{n+1}\cdots a_{n+p})^{\infty} $ such that  $$\{i_{1}=1,\,i_{2},\ldots,\,i_{t}\}\subseteq\{1,\,2,\ldots,\,n\},\quad
\mbox{and  } \{j_{1},\,j_{2},\ldots,\,j_{u}\}\subseteq\{n+1,\,\ldots,\,n+p\}$$ where $a_{i_{v}}=a_{j_{w}}=1$ for $1\leq v\leq t  $,  $1\leq w \leq u$. Now we have the following.
\begin{thm}\label{zeta} 
Let $ X(S)= $ASS$(X_{\beta}) $ for some $\beta\in (1,\,2]$. If $ X_{\beta} $ is SFT, then 
\begin{equation}\label{zeta SFT} 
\zeta_{\sigma_\beta}(r)=\zeta_{\sigma_S}(r).
\end{equation}\label{zeta sofic}
If $X_{\beta} $ is not SFT, then 
\begin{equation}
\zeta_{\sigma_\beta}(r)=(1-r)\zeta_{\sigma_S}(r).
\end{equation}
 Furthermore, in the case of SFT,
\begin{equation}
\zeta_{\sigma_\beta}(r)=\frac{1}{1-r^{i_1}-r^{i_{2}}-\cdots-r^{i_{t}}},
\end{equation}
and for strictly sofic,
\begin{equation}
\zeta_{\sigma_\beta}(r)=\frac{1}{(1-r^{i_1}-r^{i_{2}}-\cdots-r^{i_{t}})(1-r^{p})-(r^{j_{1}}+\cdots+r^{j_{u}})}.
\end{equation}
\end{thm}
\begin{proof}
First let $X_\beta$ be an SFT shift. Also let $S=\{0,\,i_{2}-1,\,\ldots,\,i_{t-1}-1,\,i_{t}-1\}$; then by Theorem (\ref{finite equivalent}), $ X(S) =$ASS$(X_\beta)$. Since $X_\beta$ and $X(S)$ are conjugate, they have the same zeta function, that is 
$$\zeta_{\sigma_\beta}(r)=\frac{1}{f_S(r^{-1})}=\frac{1}{1-r^{i_1}-r^{i_{2}}-\cdots-r^{i_{t}}}$$
where $f_{S}(x)=1-\sum_{s_{n} \in S}\frac{1}{x^{s_{n}+1}}$
\cite[Theorem 2.3]{A2}. 

Now suppose $X_\beta$ is an strictly sofic shift and let $1_{\beta}=a_{1}a_{2}\cdots a_{n}(a_{n+1}\cdots a_{n+p})^{\infty}.$ An arbitrary periodic point $x \in  X_{\beta} $ has one presentation in 
$ \mathcal G_{\beta}$ unless $$ x=(a_{n+1}\cdots a_{n+p})^{\infty} $$ where then it has exactly two
 presentations. This fact can be
 deduced from the proof of \cite[Proposition 4.7]{J}. So if $ m=pk $ $ (k \in \mathbb N) $, then every point in $ X_{\mathcal G} $ of period $ m $ is the image of exactly one
 point in $ X_{G} $ of the same period, except $ p $ points in the cycle of $ (a_{n+1}\cdots a_{n+p})^{\infty}  $ which are the image of two points of period $ m $. 
As a result, $p_{m}(\sigma_{\mathcal G_{\beta}})=  p_{m}(\sigma_{G_{\beta}})-p$ where $ p_{m}=|P_{m}| $ and $ P_{m} $ is the set of periodic points of period $ m $. When $p$ does not divide $m$, 
$ p_{m}(\sigma_{\mathcal G_{\beta}})=  p_{m}(\sigma_{G_{\beta}})$.
 Therefore,
\begin{eqnarray*}
\zeta_{\sigma_{\beta}}(r)&=&\exp\left( \sum_{m=1\atop{p\not| m}}^{\infty}\frac{ p_{m}(\sigma_{ G_{\beta}})}{m}r^{m}+\sum_{m=1\atop{p| m}}^{\infty}\frac{ p_{m}(\sigma_{G_{\beta}})-p}{m}r^{m}\right)\\
&=&\exp\left( \sum_{m=1}^{\infty}\frac{ p_{m}(\sigma_{ G_{\beta}})}{m}r^{m}-p\sum_{m=1\atop{p| m}}^{\infty}\frac{r^{m} }{m}\right)\\
&=&\zeta_{\sigma_{G_{\beta}}}(r)\times(1-r^{p}).
\end{eqnarray*}
But $ G_{\beta}\cong G_{S} $ for $ S $ as in \eqref{S}. Therefore by \cite[Theorem 2.3]{A2},
\begin{eqnarray*}
\zeta_{\sigma_{G_{\beta}}}(r)&=&\frac{1}{(1-r^{p})f_S(r^{-1})}\\
&=&\frac{1}{(1-r^{i_1}-r^{i_{2}}-\cdots-r^{i_{t}})(1-r^{p})-(r^{j_{1}}+\cdots+r^{j_{u}})}.
\end{eqnarray*}
It remains to consider the case when $ X_{\beta} $ is not sofic. We claim that $ P_{n}(X_{S})=P_{n}(X_{\beta}) +1  $ for all $ n \in \mathbb N $. 

Observe that one may assume that the initial vertex of $ \pi $, a cycle in the graph of $ G_{\beta} $, is $ \alpha_{1} $ as in Figure 4.1. Now let $ x=v^{\infty} \in P_{n}(X_{\beta}) $ with $ v=v_{1}\cdots v_{n} \in \mathcal B_{n}(X_{\beta}) $. Pick $ \pi_{\beta} $ a cycle in $ G_{\beta} $ such that $ v=\mathcal L_{\beta}(\pi_{\beta}) $ and set $ \pi_{S} $ to be the associated cycle to $ \pi_{\beta} $ in $ G_{S} $ and also let $ w=\mathcal L_{S}(\pi_{S})$. Then $ w^{\infty} \in P_{n}(X(S))$. Now define $ \varphi_{n}: P_{n}(X_{\beta})\setminus P_{1}(X_{\beta})\rightarrow P_{n}(X(S))\setminus P_{1}(X(S))  $ for all $ n\geq 2 $ such that $ \varphi_{n}(v^{\infty})=w^{\infty} $. Clearly, $ \varphi_{n} $ is well-defined. Also it is one-to-one; otherwise, for $ w^{\infty} \in P_{n}(X(S))  $, there are two different cycles $ \pi_{S} $ and $ \gamma_{S} $ such that $ w=\mathcal L_{S}(\pi_{S})=\mathcal L_{S}(\gamma_{S}) $. Let $ a_{i} $ and $ a_{j} $ be the rightmost occurance of the edges that represent the coefficients of $ 1_{\beta}$ in $ \pi_{S} $ and $ \gamma_{S} $ respectively. Then adopting notation in Figure 4.1, the terminal vertex of the edge representing $ a_{i} $ (resp. $ a_{j} $) is $ \alpha_{i+1} $ (resp. $ \alpha_{j+1} $). 

By Fischer cover of $ X(S) $, $ \alpha_{i+1}\neq \alpha_{j+1} $. Suppose $ i<j $ and $ \mathcal L_{\infty}= {\mathcal L_{S}}_{\infty } $. Then $ \mathcal L_{\infty}(\pi_{S}^{\infty})=\mathcal L_{\infty}(\gamma_{S}^{\infty}) $ implies that $ a_{i+1}=1 $. On the other hand, there is another edge starting at $\alpha_{i+1} $ and terminating at $ \alpha_{1} $ whose label is $ 1 $. Since the Fischer cover is right-resolving, it is a contradiction and so $ \varphi_{n} $ is one-to-one. But $ P_{1}(X_{\beta})=\{0^{\infty}\}  $ while $ P_{1}(X(S))=\{0^{\infty}, 1^{\infty}\}  $. So the claim is proved and we have
\begin{eqnarray*}
\zeta_{\sigma_{S}}(r)&=&\exp\left( \sum_{m=1}^{\infty}\frac{ p_{m}(\sigma_{G_{S}})}{m}r^{m}\right)\\
&=&\exp\left( \sum_{m=1}^{\infty}\frac{ p_{m}(\sigma_{G_{\beta}})+1}{m}r^{m}\right)=\zeta_{\sigma_{G_{\beta}}}(r)\times\frac{1}{1-r}.
\end{eqnarray*}
\end{proof}
Now we investigate the frequency of associated $ S $-gap shifts in the space of all $ S $-gap shifts by using topology of $ S $-gap shifts given in \cite{A1}. This topology is obtained by assigning   a real number
$x_{S}=[d_{0};\,d_{1},\,d_{2},...]$, where $[d_{0};\,d_{1},\,d_{2},...]$ is the continued fraction expansion of $x_S$, to any  $X(S)$ with $d_{0}=s_{0}$ and
$d_{n}=s_{n}-s_{n-1}$. 
 By that, a one-to-one correspondence between the $S$-gap  shifts up to conjugacy and
${\mathcal R}={\mathbb R}^{\geq 0} \setminus \{\frac{1}{n}: n \in \mathbb N \}$, up to homeomorphism, will be established and the subspace topoloy of $\mathcal R$ together with its measure structure will be induced on the space of all $S$-gap shifts.  
\begin{thm}\label{cantor set} 
Let $ \mathcal S $ be the set of all  $ S $-gap shifts associated to some $ X_{\beta} $. Then
 $ \mathcal S $ is a Cantor set on the space of all $ S $-gap shifts (a nowhere dense perfect set).
Entropy is a complete invariant for the conjugacy classes of $\mathcal{S}$.
\end{thm}
\begin{proof}
First suppose $ X(S) $ does not satisfy \eqref{varepsilon} and $ x_{S}=[d_{0};\,d_{1},\,\ldots] $ corresponds to $ X(S) $ \cite{A1}. Let $ N  $ be the least integer such that $d_{N}d_{N+1}\cdots< d_{1}d_{2}\cdots$ and set $\gamma_{i}:=[d_{0};d_{1},\ldots,d_{i}] $, $ i \in \mathbb N_{0} $. If $ N $ is even, set $ U=(\gamma_{N}, \gamma_{N+1})$ and otherwise, $ U=(\gamma_{N+1}, \gamma_{N})$. Then no points of $ U $ satisfies \eqref{varepsilon} and so none is an associated $ S $-gap shift. This shows $ \mathcal S $ is closed. 

Now let $ X(S) \in \mathcal S $ and $ V $  be a neighborhood of $ x_{S} $. Note that two real numbers are close if sufficiently large numbers of their first partial quotients in their  continued fraction expansion are equal. So we can select two points $ x_{S'},
 x_{S''} \in V$ such that $ X(S') $ satisfies \eqref{varepsilon} and $ X(S'') $ does not satisfy \eqref{varepsilon}. This implies that all points of $ \mathcal S $ are limit points of themselves and $ \mathcal S $ is nowhere dense.

The second part follows from the fact that the 
 entropy is a complete invariant for the conjugacy classes of $ \beta $-shifts.
\end{proof}

The interesting equivalencies happening in case of sofics is when the two systems under investigation have a common SFT extension. This clearly cannot be done for non-sofics. The most natural extension of this idea is when two non-sofic systems have a common synchronized extension and in particular, when the legs are right-resolving. This has been studied in \cite{F'} and \cite{F}. Recall that  when $ X $ is a synchronized system with Fischer cover $ \mathcal G=(G, \mathcal L) $, then $ \mathcal L_{\infty} $ is a.e. 1-1 \cite{FF}. So for our case we have the following.
\begin{thm}\label{common synchronized} 
$ X_{\beta} $ and ASS$(X_{\beta}) $ have a common synchronized 1-1 a.e. extension with right-resolving legs.
\end{thm}

\bibliographystyle{amsplain}

\end{document}